\documentclass[reqno,11pt]{amsart}
\usepackage[english]{babel}
\usepackage{amsmath}
\usepackage{amssymb}
\usepackage{enumerate}
\usepackage{a4wide}
\usepackage{hyperref}
\usepackage{color}
\usepackage{graphicx}
\usepackage{subfigure}

\newtheorem{thm}{Theorem}[section]
\newtheorem{prop}[thm]{Proposition}

\newtheorem{lemma}[thm]{Lemma}
\newtheorem{cor}[thm]{Corollary}
\newtheorem{remark}[thm]{Remark}

\numberwithin{equation}{section}
\renewcommand{\proof}{\smallskip\noindent{\it Proof. }}
\renewcommand{\endproof}{\qed\smallskip}
\newcommand{\R}{\mathbb{R}}

\newcommand{\Rd}{\mathbb{R}^{d}}

\newcommand{\pp}{\mathcal{P}}

\newcommand{\e}{\varepsilon}

\newcommand{\rhot}{\widetilde{\rho}}
\newcommand{\rhob}{\overline{\rho}}

\let\oldmarginpar\marginpar
\renewcommand\marginpar[1]{\-\oldmarginpar[\footnotesize #1]{\footnotesize #1}}

\DeclareMathOperator{\dive}{div}

\newcommand{\dep}{\partial}

\begin{document}

\title[Quadratic diffusion equations with long-range attraction]{Stationary states of quadratic diffusion equations with long-range attraction}

\author{M. Burger}
\address{Martin Burger, Institute for Computational and Applied Mathematics, Westf\"alische Wilhelms Universit\"at (WWU) M\"unster, Einsteinstraße 62, D 48149 M\"unster, Germany}
\email{martin.burger@wwu.de}
\author{M. Di Francesco}
\address{Marco Di Francesco, Department of Pure and Applied Mathematics, University of L'Aquila, Via Vetoio Loc. Coppito, 67100 L'Aquila, Italy}
\email{mdifrance@gmail.com}
\author{M. Franek}
\address{Marzena Franek, Institute for Computational and Applied Mathematics, Westf\"alische Wilhelms Universit\"at (WWU) M\"unster, Einsteinstraße 62, D 48149 M\"unster, Germany}
\email{marzena.franek@uni-muenster.de}

\date{\today}

\begin{abstract}
We study the existence and uniqueness of nontrivial stationary solutions to a nonlocal aggregation equation with quadratic diffusion arising in many contexts in population dynamics. The equation is the Wasserstein gradient flow generated by the energy $E$, which is the sum of a quadratic free energy and the interaction energy. The interaction kernel is taken radial and attractive, nonnegative and integrable, with further technical smoothness assumptions. The existence vs. nonexistence of such solutions is ruled by a threshold phenomenon, namely nontrivial steady states exist if and only if the diffusivity constant is strictly smaller than the total mass of the interaction kernel. In the one dimensional case we prove that steady states are unique up to translations and mass constraint. The strategy is based on a strong version of the Krein-Rutman theorem. The steady states are symmetric with respect to their center of mass $x_0$, compactly supported on sets of the form $[x_0 -L, x_0+L]$, $C^2$ on their support, strictly decreasing on $(x_0,x_0+L)$. Moreover, they are global minimizers of the energy functional $E$. The results are complemented by numerical simulations.
\end{abstract}

\maketitle


\section{Introduction}

Phenomena with long-range aggregation and short-range repulsion arise in many instances in population biology such as chemotaxis of cells, swarming or flocking of animals. A variety of mathematical models has been proposed for such situations, at the particle as well as at the continuum (mean field) level. In particular, if nonlocal \emph{repulsion} acts at a smaller scale with respect to nonlocal \emph{attractive} forces in the large particle limit, then a nonlocal repulsion term can be replaced by a \emph{local} term with nonlinear diffusion, we refer to \cite{mogilner_edelstein,boi_capasso_morale,topaz_bertozzi,topaz_bertozzi_lewis,Chuang_dorsogna_marthaler_bertozzi_chayes,
dorsogna_chuang_bertozzi_chayes,kang_perthame_stevens,primi_stevens_velazquez,bodnar_velazquez,li_zhang} for several examples. A prototype model, which we shall also investigate further in this paper, is given by
\begin{equation}\label{eq:main_intro}
    \dep_t \rho = \dive \left(\rho \nabla \left( \e \rho - G*\rho \right)\right)
\end{equation}
where the convolution is carried out with an aggregation kernel $G$ such that $G(x)=g(|x|)$ with $g'(r)>0$ as $r>0$. This model arises in a natural way as the limit of a stochastic interacting particle model with pair interactions (cf. \cite{morale_capasso_oelschlaeger,oelschlaeger} respectively \cite{golse} for general background). Models with the same structure have been recently used to model opinion formation, cf. \cite{porfiri,sznajd}. Here we shall restrict to the case in which \eqref{eq:main_intro} is posed on the whole space $\R^d$.

In case $G$ is $\lambda$-convex, then the equation \eqref{eq:main_intro} can be formulated as a gradient flow in the Wasserstein metric (cf. \cite{AGS,Villani_book,villani_book_2}) of the energy (or entropy) functional
\begin{equation}
    E[\rho]:=\frac{\varepsilon}2\int_{\R^d} \rho^2(x) dx - \frac{1}{2}\int_{\R^d}\int_{\R^d} G(x-y) \rho(y) \rho(x) dy dx, \label{entropyfunctional}
\end{equation}
see \cite{mccann,BenedettoCagliotiPulvirenti97,carrillo_mccann_villani1,carrillo_mccann_villani2}. Models of the above form have been investigated with respect to several aspects, e.g. existence and uniqueness in the context of entropy solutions \cite{burger_capasso_morale,bertozzi_slepcev}, well-posedness in the context of Wasserstein gradient flows as a special case of the theory developed in \cite{AGS}, and - in particular in connection with the classical Patlak-Keller-Segel model for chemotaxis \cite{patlak,keller_segel} - with respect to blow-up vs. large-time existence, cf. e. g. \cite{jaeger_luckhaus,herrero_velazquez,blanchet_dolbeault_perthame,corrias_perthame_zaag} for models with linear diffusion and \cite{kowalczyk,Calvez_Carrillo} for models with nonlinear diffusion.

An interesting and important question is the characterization of large-time behavior of solutions to equations of the form \eqref{eq:main_intro}, which is related with the possible existence of nontrivial steady states, even when the quadratic diffusion is replaced by a more general nonlinear diffusion. This issue is solved in detail for purely diffusive equations, in which solutions decay to zero with a prescribed rate for large times and behave like the (compactly supported) Barenblatt profiles, cf. the classical works of Vazquez on the self-similar behavior of the porous medium equation, which are nicely collected in the book \cite{Vazquez}, as well as the papers \cite{otto,carrillo_toscani}. In the purely nonlocal case, namely when $\e=0$, this issue has been studied extensively in many papers \cite{LiTos,BuDiF,laurent,bertozzi_laurent,bertozzi_carrillo_laurent,Carrillo_DiFrancesco_Figalli_et_al,bertozzi_brandman,bertozzi_laurent_rosado,
huang_bertozzi,fellner_raoul1,fellner_raoul2}, combined with the study of the regularity of solutions compared to the \emph{attractive singularity} of the interaction kernel. In particular, solutions are known to concentrate to a Dirac delta centered at the initial center of mass (invariant) either in a finite or in an infinite time, depending on the properties of the kernel $G$ at $x=0$. When the kernel $G$ is supported on the whole space, the Dirac delta is the unique steady state of \eqref{eq:main_intro} with unit mass and zero center of mass.

The asymptotic behavior in the general case with both nonlinear diffusion and nonlocal interaction has been only partially addressed. A first attempt in this direction was performed in \cite{BuDiF}, in which the existence of steady states of \eqref{eq:main_intro} for sufficiently small $\e$ and the non-existence for large $\e$ in the one-dimensional case was proven by means of the pseudo-inverse representation of the Wasserstein distance. More refined results in a similar model derived in \cite{Hillen_Painter,Painter_Hillen} with cut-off density have been found in \cite{Burger_DiFrancesco_DolakStruss,Laurencot_Wrozsek}.
Parallel to their work, the authors of the present paper recovered the results in \cite{bedrossian}, in which a quasi sharp result of existence of minimizers for the energy $E[\rho]$ in a multi dimensional framework has been proven.

A key open (to our knowledge) problem in this context is the \emph{uniqueness} of steady states under mass and center of mass constraint, its main difficulty being the fact that the functional $E$ is neither convex in the classical sense nor in the displacement convex sense \cite{mccann} (except when $G$ is concave on $\R^d$, see \cite{carrillo_mccann_villani1}) when $\e<\int G$ with $G\geq 0$.

In this paper we further investigate the detailed structure of steady states in one space dimension. We remark that our work does not go into the direction of providing sharp regularity conditions on the kernel $G$. Roughly speaking, $G$ is smooth, radial (with decreasing profile), nonnegative, integrable and supported on the whole space. The precise assumptions on $G$ are stated at the beginning of Section \ref{sec:preliminaries}. We found out that for a kernel $G$ decreasing on $x>0$, the $L^1$-norm of $G$ compared to $\e$ marks a threshold:
\begin{itemize}
\item If $\e \geq \int G $, then there exists no steady state.
\item If $\e < \int G$, there exists a stationary state.
\end{itemize}
The same results are recovered in \cite{bedrossian}, except for the critical case $\e=\int G$.

The main result of our paper deals with the case $d=1$. Here we can give a detailed characterization of the stationary states in the case $\e < \int G$.
If $G'$ only vanishes at zero, then there exists a \emph{unique} stationary state (up to translation), which is a minimizer of the energy $E$ at fixed mass. The stationary states have compact support, which increases with $\e$. Moreover, such a steady state is symmetric with a single maximum at the center of mass, and monotone on both sides of the center of mass. The precise statement can be found in Theorem \ref{thm:main1d}. If $\e$ is small enough, the stationary states are concave on their support. The main tool in the proof of the main result is the statement of the stationary equation as an eigenvalue problem, see formula \eqref{eq:stat_compact}, together with the use of a strong version of Krein-Rutman theorem (cf. Theorem \ref{thm:KR2}), which allows to characterize the steady states as eigenfunctions corresponding to a \emph{simple} eigenvalue.

The uniqueness of the stationary state is somewhat surprising, since the energy functional $E$ is not convex for $\e < \int G$ and thus one might expect other stationary points of $E$. On the other hand also in the case $\e = 0$ one can see that there is a unique (measure) stationary state concentrated at the center of mass if $G$ has global support. Adding the squared norm for positive $\e$ makes the functional closer to convex, and thus probably does not lead to additional stationary points.

The paper is organized as follows. In Section \ref{sec:preliminaries} we recall the statement of the problem and provide some preliminary regularity results. In Section \ref{sec:stationary_multi_d} we complement our results with those proven in \cite{bedrossian} and provide sharp conditions on $\e$ and $\int G$ for the existence of non trivial steady states. In Section \ref{sec:stationary_one_d} we prove our main results about the uniqueness of steady states in one space dimension. Finally, in Section \ref{sec:numerics} we complement our results with some numerical simulations.

\section{Preliminaries}\label{sec:preliminaries}

We consider the evolution equation
\begin{equation}\label{eq:main_evo}
    \dep_t \rho = \dive \left(\rho \nabla \left( \e \rho - G*\rho \right)\right)
\end{equation}
and its stationary version
\begin{equation}\label{eq:main_stat}
    0= \dive \left(\rho \nabla \left( \e \rho - G*\rho \right)\right)
\end{equation}
posed on the whole space $\Rd$. Due to the applied setting of the model, we shall consider here only nonnegative solutions $\rho \geq 0$.

\textbf{Assumptions on the kernel $G$}. We shall assume throughout the paper that the interaction kernel $G$ satisfies
\begin{enumerate}
  \item $G\geq 0$, and $\mathrm{supp}(G)=\R^d$,
  \item $G\in W^{1,1}(\Rd)\cap L^\infty(\Rd)\cap C^2(\Rd)$,
  \item $G(x)=g(|x|)$ for all $x\in \Rd$,
  \item $g'(r)<0$ for all $r>0$,
  \item $g''(0)<0$,
  \item $\lim_{r\rightarrow +\infty} g(r) =0$.
\end{enumerate}

We emphasize here that providing sharp conditions on the regularity of $G$ is not a purpose of the present paper.  Let us recall that the equation \eqref{eq:main_evo} preserves the total mass
\begin{equation*}
    M=\int \rho (x) dx
\end{equation*}
and center of mass
\begin{equation*}
    CM[\rho(t)] : = \int x \rho (x,t) dx.
\end{equation*}
Moreover, it is easily seen that, for a given stationary state $\rho$ solving \eqref{eq:main_stat}, $M\rho$ and $\rho(\cdot +x_0)$ are still stationary states for all $M>0$ and $x_0\in \R$.

We shall therefore assume $M=1$ for simplicity. For future use we introduce the space
$$\mathcal{P}=\left\{\rho \in L^1_+(\Rd):\ \int_{\Rd}\rho(x) dx = 1\right\}. $$

Moreover, from now on we shall assume for simplicity
\begin{equation*}
    \|G\|_{L^1}=\int G(x) dx =1.
\end{equation*}
This is not restrictive since the kernel $G$ can always be normalized by modifying the diffusion constant and the time scale as follows
\begin{equation*}
    \partial_\tau \rho = \mathrm{div}\left(\rho\nabla(\e' \rho - \widetilde{G}*\rho)\right),\quad \tau = \|G\|_{L^1(\R^d)},\quad \widetilde{G} = G/\|G\|_{L^1(\R^d)},\quad \e' = \e/\|G\|_{L^1(\R^d)}.
\end{equation*}

Let us recall the following results on the existence and uniqueness of gradient flow solutions to \eqref{eq:main_evo}, which follows from the theory developed in \cite{AGS}. In this sense, let us introduce one of the most important tools related with the study of the evolution equation \eqref{eq:main_evo} and in particular with the stationary version \eqref{eq:main_stat}, namely the \emph{energy functional}
\begin{equation*}
    E[\rho]:=\frac{\varepsilon}{2}\int_{\R^d} \rho^2(x) dx - \frac{1}{2}\int_{\R^d}\int_{\R^d} G(x-y) \rho(y) \rho(x) dy dx.
\end{equation*}

\begin{thm}[\cite{AGS}]\label{thm:ags}
Let $\rho_0\in L^2\cap\mathcal{P}$ such that $\rho\geq 0$ and $E[\rho_0]<+\infty$. Let $G$ satisfy the above assumptions. Then there exists a unique weak solution $\rho$ to \eqref{eq:main_evo} with
\begin{itemize}
  \item $E[\rho(t)]<+\infty$ for all $t\geq 0$.
  \item $\sqrt{\rho}\nabla(\varepsilon\rho - G*\rho) \in L^2([0,T]\times \R^2)$ for all $T>0$
\end{itemize}
such that the following energy identity is satisfied
\begin{equation}\label{eq:energy_identity}
    E[\rho(t)] +\int_0^T \int_{\R^d}\rho\left|\nabla(\varepsilon \rho - G*\rho)\right|^2 dx dt = E[\rho_0].
\end{equation}
\end{thm}

In particular, the equation \eqref{eq:main_evo} produces the following regularizing effect.
\begin{lemma}[Regularity of $L^2\cap \mathcal{P}$ steady states]\label{lem:reg}
Let $\rho_0\in L^2\cap\mathcal{P}$. Then, the corresponding solution $\rho(t)$ to \eqref{eq:main_evo} satisfies
\begin{equation}\label{regularity}
    \int \rho \left|\nabla \rho\right|^2 dx < +\infty,
\end{equation}
for almost every $t>0$. In particular, let $\rho$ be an $L^2\cap\mathcal{P}$ steady state to \eqref{eq:main_evo}, then $\rho$ satisfies \eqref{regularity} and $\rho \in C^2$ on $\mathrm{supp}[\rho]$.
\end{lemma}

\proof
Due to \eqref{eq:energy_identity}, the quantity
\begin{equation*}
    \varepsilon \int \rho \left|\nabla \rho\right|^2 dx - 2\varepsilon\int \rho \nabla \rho \cdot \nabla G *\rho dx + \int \rho \left|\nabla G*\rho\right|^2 dx
\end{equation*}
is finite for almost every $t>0$, and therefore, in view of Cauchy-Schwarz inequality, we have
\begin{equation*}
   \frac{\varepsilon}{2} \int \rho \left|\nabla \rho\right|^2 dx - C(\varepsilon) \int \rho \left|\nabla G*\rho\right|^2 dx + \int \rho \left|\nabla G*\rho\right|^2 dx < +\infty,
\end{equation*}
and thanks to the smoothness assumptions on $G$ we have the assertion \eqref{regularity}. Let $\rho$ be a steady state, then $\rho$ satisfies \eqref{regularity} too. This implies in particular that $\nabla \rho$ is almost everywhere finite on $\R^2$. The energy identity \eqref{eq:energy_identity} implies then
\begin{equation*}
    \rho\left|\nabla (\varepsilon \rho - G*\rho)\right|^2 =0
\end{equation*}
for almost every $x\in \R^2$. This means that
\begin{equation*}
    \varepsilon \rho - G*\rho = \hbox{constant}
\end{equation*}
almost everywhere on every connected component of the support of $\rho$. By convolution with standard mollifiers, one can easily see that $\e\rho-G*\rho=C$ for a given $C$ depending on the connected component of $\mathrm{supp}[\rho]$. Since $G$ is $C^2$, this easily implies $\rho \in C^2$ on $\mathrm{supp}[\rho]$.
\endproof

\begin{cor}[$1d$ regularity]\label{cor:1dreg}
Let $\rho$ be an $L^2\cap \mathcal{P}$ solution to \eqref{eq:main_stat} in one space dimension. Then $\rho$ is continuous on $\R$.
\end{cor}

\proof
Apply the result in Lemma \ref{lem:reg} to the case $d=1$. Since $\frac{d}{dx}\rho^{3/2} \in L^2$, the one dimensional Sobolev embedding implies that $\rho^{3/2}$ is continuous.
\endproof

\section{Stationary solutions in multiple dimensions}\label{sec:stationary_multi_d}

In this section we state the necessary and sufficient conditions on $\e$ and $\|G\|_{L^1}$ such that there exists non trivial steady states
\begin{equation}\label{eq:main_stat_multid}
    \rho \nabla (\e \rho - G*\rho ) =0
\end{equation}
in the set $L^2\cap \mathcal{P}$. During our work, we realized that J.Bedrossian has obtained similar results in \cite{bedrossian}, based on ideas and strategies developed in \cite{lieb} and \cite{lions}. In order to simplify the coverage of the paper, we shall state the result in \cite[Theorem 1]{bedrossian} and prove all other results. Notice that the critical case $\e=\|G\|_{L^1}$ was not covered in \cite{bedrossian}.

Let us start by focusing on the interplay between the solutions to \eqref{eq:main_stat_multid} and the variational calculus on the energy functional
\begin{equation*}
    E[\rho]:=\varepsilon\int_{\R^d} \rho^2(x) dx - \frac{1}{2}\int_{\R^d}\int_{\R^d} G(x-y) \rho(y) \rho(x) dy dx.
\end{equation*}

In the next proposition we prove that being a minimum for the entropy functional is a sufficient condition for being a solution to \eqref{eq:main_stat_multid}.

\begin{prop}[Stationary solutions via energy minimization]\label{prop:min_implies_stat}
Let $\rho\in L^2(\R^d)$ be a minimizer for the energy functional
\begin{equation*}
    E[\rho]:=\frac{1}{2}\int_{\Rd}\rho\left(\e\rho - G*\rho\right) dx
\end{equation*}
on $\mathcal{P}$.
Then
\begin{equation*}
    \rho\nabla\left(\e \rho - G*\rho\right)=0\qquad \hbox{a. e. in}\ \Rd.
\end{equation*}
\end{prop}

\proof
Let $V\in C^1_c(\Rd)$ be an arbitrary vector field and let $u(x,s)$ be a local solution to the continuity equation
\begin{equation*}
    \dep_s u(x,s) + \dive (u(x,s) V(x)) = 0
\end{equation*}
with initial datum
\begin{equation*}
    u(x,0) = \rho(x)
\end{equation*}
with $\rho$ being the minimizer for $E$ given in the hypothesis. Such a $u$ can be constructed by solving the characteristic ODE
\begin{equation*}
    \frac{d}{ds} X(x,s) = V(X(x,s))
\end{equation*}
coupled with the initial datum
\begin{equation*}
    X(x,0)=x
\end{equation*}
locally in $s=0$, with the local solution $X(x,s)$ being $C^1$, and by taking $u(x,s):=[(X(\cdot,s))_\sharp \rho] (x,s)$, i. e. $u(x,s)$ defined via
\begin{equation*}
   \int \phi(x) u(x,s)dx = \int \phi (X(x,s))\rho(x) dx,\qquad \hbox{for all}\ \ \phi\in C^1_c (\Rd)
\end{equation*}
(cf. for instance \cite[Chapter 8, Lemma 5.5.3]{AGS}). For all $s$ in the interval of existence of $u$ we have
\begin{equation*}
    \int u(x,s) dx = 1,\qquad u(x,s)\geq 0 \ \ \hbox{a.e.}
\end{equation*}
and therefore the map $s\mapsto E[u(\cdot,s)]$ has a local minimum at $s=0$. Hence
\begin{align*}
    & 0 \leq \frac{d}{ds} E[u(\cdot,s)]|_{s=0} = \int (\e u - G*u ) \dep_s u dx|_{s=0}\\
    & \ = - \int (\e u - G*u ) \dive (uV) dx|_{s=0}= \int \rho \nabla (\e \rho - G*\rho ) \cdot V dx
\end{align*}
and replacing $V$ with $-V$ we obtain
\begin{equation*}
    0 \geq \int \rho \nabla (\e \rho - G*\rho ) \cdot V dx
\end{equation*}
and therefore
\begin{equation*}
    \int \rho \nabla (\e \rho - G*\rho ) \cdot V dx = 0,\qquad \hbox{for an arbitrary}\ \ V\in C^1_c(\Rd)
\end{equation*}
which is the desired assertion.
\endproof

Let us now compute the first and the second order Gateaux derivatives of $E$.

\begin{lemma}\label{lem:derivativE}
Let $\rho\in L^2\cap \mathcal{P}$ be a solution to \eqref{eq:main_stat_multid}. Then, $\rho$ is a stationary point for the energy functional $E$. Moreover, the second order Gateaux derivative of $E$ on $\rho$ satisfies
\begin{equation}\label{eq:second_derivative}
    \frac{d^2}{d \delta^2} E[\rho + \delta v]|_{\delta=0} = \varepsilon \int_{\R^d} v^2(x) dx - \int v(x) G*v(x) dx,
\end{equation}
for all $v=\mathrm{div}(\rho V)$ and $V\in C^1_c(\R^d)$.
\end{lemma}

\proof

Suppose $\rho\in L^2\cap\mathcal{P}$ satisfies \eqref{eq:main_stat_multid}. Let us compute
\begin{equation*}
    \lim_{\delta\rightarrow 0}\frac{1}{\delta}\left( E[\rho + \delta v] - E[\rho]\right)
\end{equation*}
with $v=\dive (\rho V)$ for an arbitrary vector field $V\in C^1_c$, which implies $\int v(x) dx = 0$. We obtain
\begin{align*}
    & \frac{1}{\delta}\left( E[\rho + \delta v] - E[\rho]\right)\\
    & = \frac{\e}{2\delta}\int_{\mathrm{supp}(\rho+\delta v)}(\rho + \delta v)^2 dx - \frac{\e}{2\delta}\int_{\mathrm{supp}(\rho)}\rho ^2 dx \\
    & -\frac{1}{2\delta}\int_{\mathrm{supp}(\rho+\delta v)}(\rho+\delta v) G*(\rho + \delta v) dx + \frac{1}{2\delta}\int_{\mathrm{supp}(\rho)}\rho G*\rho dx.
\end{align*}
Therefore we easily get
\begin{align*}
   & \lim_{\delta\rightarrow 0}\frac{1}{\delta}\left( E[\rho + \delta v] - E[\rho]\right) = \int v(\e \rho - G*\rho) dx\\
   & \ = \int \dive (\rho V) (\e \rho - G*\rho) dx = - \int \rho V \cdot \nabla  (\e \rho - G*\rho) dx.
\end{align*}
Therefore, $\rho$ is a stationary point for $E$ under the constraint $\int \rho dx =1$. The computation of the second derivative of $E$ on $\rho$ yields
\begin{align*}
    & \frac{d^2}{d\delta^2} E[\rho + \delta v] = \frac{d^2}{d\delta^2} \frac{\e}{2}\int_{\mathrm{supp}(\rho+\delta v)}(\rho + \delta v)^2 dx\\
    & \ \ -\frac{d^2}{d\delta^2}\frac{1}{2}\int_{\mathrm{supp}(\rho+\delta v)}(\rho+\delta v) G*(\rho + \delta v) dx\\
    & \ =  \e\int v^2 dx -\int v  G*v dx
\end{align*}
which is independent on $\delta$ and therefore it is valid also on $\delta=0$.
\endproof

Before we start analyzing the existence or nonexistence of steady states, we introduce a very simple technical lemma which will be very useful in the sequel.

\begin{lemma}\label{lem:trick}
Suppose $\rho\in L^2\cap\mathcal{P}$ is a solution to \eqref{eq:main_stat_multid} having connected support. Then
\begin{equation*}
    \e\rho(x) = \int_{\mathrm{supp}[\rho]} G(x-y)\rho(y) dy + C
\end{equation*}
for all $x\in \mathrm{supp}[\rho]$ with $C=2E[\rho]$. Moreover, in case $\mathrm{supp}[\rho]$ has infinite measure, then $C=E[\rho]=0$.
\end{lemma}

\proof
It is immediate from \eqref{eq:main_stat_multid} that
\begin{equation}\label{eq:support}
    \e\rho(x) = \int_{\mathrm{supp}[\rho]} G(x-y)\rho(y) dy + C
\end{equation}
for all $x\in \mathrm{supp}[\rho]$ for a certain constant $C$. Then, we multiply \eqref{eq:support} by $\rho(x)$ and integrate over $\mathrm{supp}[\rho]$ to obtain
\begin{equation*}
    \e\int_{\mathrm{supp}[\rho]}\rho^2(x) dx = \int_{\mathrm{supp}[\rho]}\int_{\mathrm{supp}[\rho]} G(x-y)\rho(y)\rho(x) dy dx + C,
\end{equation*}
where we have used that $\rho$ has unit mass. It is therefore clear that $C= 2 E[\rho]$. Suppose now that $\mathrm{supp}[\rho]$ has infinite measure. Suppose by contradiction that $C\neq 0$. Let $\{x_k\}\subset \mathrm{supp}[\rho]$ be a sequence of points such that $|x_k|\rightarrow +\infty$. We have, for all $k$,
\begin{equation*}
    \e \rho(x_k) - \int_{\mathrm{supp}[\rho]} G(x_k - y) \rho(y) dy = C
\end{equation*}
and therefore the same expression should hold in the limit $k\rightarrow +\infty$. Now, the assumptions on $G$ imply that the integral
\begin{equation*}
    \int_{\mathrm{supp}[\rho]} G(x_k - y) \rho(y) dy
\end{equation*}
converges to zero as $k\rightarrow +\infty$. This is due to Lebesgue's dominated convergence Theorem. Therefore, the term $\rho(x_k)$ has a limit $C$ as $k\rightarrow +\infty$. Such limit is the same for all diverging sequences of points $\{x_k\}\subset \mathrm{supp}[\rho]$, which means
\begin{equation*}
    \lim_{x\in\mathrm{supp}[\rho], |x|\rightarrow +\infty} \rho( x) = C.
\end{equation*}
Now, since $\mathrm{supp}[\rho]$ has infinite measure, then $C\neq 0$ implies that $\rho$ is not integrable, which is a contradiction. Therefore $C=0$.
\endproof

\subsection{Non existence of nontrivial steady states for $\e>1$}

We start by covering the case $\e>1$. Here, there exist no nontrivial steady states, as it follows from the following simple lemma.

\begin{lemma}\label{lem:supercritical}
Let $\e>1$. Then, there exists no stationary solutions to \eqref{eq:main_stat_multid} in the space $L^2\cap \mathcal{P}$.
\end{lemma}

\proof
We first prove that there exists no minimizer for $E[\rho]$ under the mass constraint $\int \rho = 1$ and $\rho\geq 0$. To see this, we use Young inequality for convolutions as follows
\begin{equation}\label{est:convexcase}
    E[\rho] = \frac{\e}{2}\int \rho^2 dx - \frac{1}{2}\int \rho G*\rho dx \geq \frac{\e}{2}\int \rho^2 dx - \frac{\|G\|_{L^1}}{2}\int \rho^2 dx = \frac{\e- 1}{2}\int \rho^2 dx
\end{equation}
with $\e-1>0$. Moreover, we have the simple estimate $E[\rho]\leq C\|\rho\|_{L^2}^2$. Take a family of functions $\rho_\lambda(x)\geq 0$ such that $\int \rho_\lambda(x) dx =1$ and $\int \rho_\lambda^2 (x) dx\rightarrow 0$ as $\lambda \rightarrow +\infty$. To construct such a family, we just take a fixed $L^2_+(\Rd)$ function $\rho\not\equiv 0$ and rescale it by $\rho_\lambda(x)=\lambda^{-d}\rho(\lambda^{-1} x)$. For such a family we therefore have
\begin{equation*}
    E[\rho_\lambda]\rightarrow 0,\qquad \hbox{as}\ \ \lambda\rightarrow \infty.
\end{equation*}
Therefore, it is impossible to have a minimizer $\rho_\infty$ for $E[\rho]$ in the set $\left\{\rho\in L^1_+ :\ \int \rho = 1\right\}$ because \eqref{est:convexcase} would imply that $E[\rho_\infty]>0$ and we would necessarily have $0<E[\rho_\lambda]<E[\rho_\infty]$ for $\lambda$ large enough.

Now we prove that there exist no steady states. Suppose by contradiction that $\rho$ is a steady state. Then, due to Lemma \ref{lem:derivativE} $\rho$ is a stationary point for $E$. Moreover, the formula \eqref{eq:second_derivative} implies that the functional $E$ is strongly convex, and therefore admits only one stationary point, which coincides with its global minimizer. But this contradicts the non existence of a global minimizer proven above.
\endproof

\subsection{The critical case $\e=1$}

We aim to solve
\begin{equation}\label{eq:main_stat_multid_1}
    0= \dive \left(\rho \nabla \left(\rho - G*\rho \right)\right).
\end{equation}

We shall prove that no $L^2\cap\mathcal{P}$ steady states exist in this case.

\begin{thm}[Non-existence of nontrivial steady states for $\e=1$]\label{thm:nonexist1}
There exists no solutions to \eqref{eq:main_stat_multid_1} in $L^2\cap\mathcal{P}$.
\end{thm}

\proof
From Cauchy--Schwartz inequality we know that
\begin{equation*}
    \int_{\R^d}\rho G*\rho dx \leq \|\rho\|_{L^2(\R^d)} \|G*\rho\|_{L^2(\R^d)}
\end{equation*}
and the equality in the above formula holds if and only if $\rho$ and $G*\rho$ are proportional.
In terms of the functional $E$ this means that
\begin{equation*}
    E[\rho]\geq 0\qquad \hbox{for all}\ \ \rho\in L^2(\R^d)\cap\mathcal{P}(\R^d).
\end{equation*}
As in Lemma \ref{lem:supercritical}, we have the estimate $E[\rho]\leq C\|\rho\|_{L^2}^2$, and using once again the family $\rho_\lambda$ of Lemma \ref{lem:supercritical} we see that $\inf_{\rho\in L^2\cap\mathcal{P}}E[\rho] = 0$. Assume by contradiction that there exists a stationary solution $\rho_\infty$. Then, due to the result in Lemma \ref{lem:derivativE} and in view of Cauchy--Schwartz inequality, the second order derivative of $E$ is nonnegative everywhere. Hence, the functional $E$ is convex and therefore $\rho_\infty$ is a global minimizer for $E$ under the constraint $\rho\in L^2\cap \mathcal{P}$. Then, we must have $E[\rho_\infty]=0$, which means that $\rho$ and $G*\rho$ are proportional, i. e. there exists a constant $\lambda\in \R_+$ such that
\begin{equation}\label{eq:proportional}
    \rho_\infty(x) = \lambda G*\rho_\infty(x)
\end{equation}
almost everywhere on $\R^d$. Integrating \eqref{eq:proportional} over $\R^d$ yields
\begin{equation*}
    1 = \lambda \|G\|_{L^1(\R^d)}=\lambda
\end{equation*}
and hence
\begin{equation}\label{eq:stat1}
    \rho_\infty(x) = G*\rho_\infty (x)
\end{equation}
almost everywhere on $\R^d$. We can then apply the Fourier transform
\begin{equation*}
    \widehat{f}(\xi)=\int_{\R^d}e^{-2\pi i x\cdot \xi} f(x) dx
\end{equation*}
to both members of the equation \eqref{eq:stat1} to obtain
\begin{equation*}
    \widehat{\rho_\infty}(\xi)=\widehat{G}(\xi)\widehat{\rho_\infty}(\xi),\qquad \xi\in\R^d.
\end{equation*}
We have
\begin{equation*}
    |\widehat{G}(\xi)|\leq \int_{\R^d}| G(x)| dx = 1.
\end{equation*}
Moreover, since $G$ is even, then $\widehat{G}(\xi)< 1$ for all $\xi\neq 0$. In order to see that, write
\begin{equation*}
    \widehat{G}(\xi)=\int_{\R^d}\prod_{k=1}^d e^{-2\pi i x_k \xi_k} G(x) dx = \int_{\R^d}\prod_{k=1}^d (\cos(2\pi x_k \xi_k) - i\sin(2\pi x_k \xi_k) ) G(x) dx,
\end{equation*}
then $G$ being even easily implies that only real valued contributions survive in the above integral; such real valued contributions are of the form
\begin{equation*}
    \int_{\R^d} f_{h,k}(x,\xi) G(x) dx
\end{equation*}
where the functions $f_{h,k}$ are such that $|f_{h,k}(x,\xi)|\leq 1$ and $|f_{h,k}(x,\xi)|<1$ for $x$ ranging on a set of positive measure. Therefore, we have proven that
\begin{equation*}
    \widehat{\rho_\infty}(\xi)=0\quad \hbox{for all}\quad \xi\neq 0
\end{equation*}
and $\widehat{\rho_\infty}(0)=1$. This implies that $\rho(x)=0$ almost everywhere, which contradicts the fact that $\rho$ has unit mass.
\endproof

\subsection{Stationary solutions for $\e<1$}

Let us now provide a minimizer for the entropy functional in the case $\e<1$, which implies the existence of a nontrivial $L^2\cap\mathcal{P}$ steady state for \eqref{eq:main_stat_multid} in view of Proposition \ref{prop:min_implies_stat}. Such result is proven rigorously in \cite[Theorem 1]{bedrossian}, which we recall here.

\begin{thm}[Existence of minimizers, \cite{bedrossian}]\label{thm:existence_minimizers}
Let $\e<1$. Then, there exists a radially symmetric non-increasing minimizer $\rho\in\mathcal{P}\cap L^2(\Rd)$ for the entropy functional $E$ restricted to $\mathcal{P}$ with $\rho\neq 0$.
\end{thm}

We send the reader to \cite{bedrossian} for the details of the proof, which is based on a sort of subadditivity property needed to provide suitable compactness of the minimizing sequence, cf. \cite[Lemma 2]{bedrossian}. For the sake of clarity, we shall still provide the simple proof of the fact that global minima of $E$ under mass constraint are strictly negative, which forces the minimizer to be non zero.

\begin{lemma}
Let $\e<1$. Then, $\inf_{\rho\in\mathcal{P}\cap L^2(\Rd)} E[\rho]<0$.
\end{lemma}

\proof
We consider the family $\sigma_\lambda\in L^2\cap \mathcal{P}$
\begin{equation*}
    \sigma_\lambda(x)=\frac{1}{2\lambda} \chi_{[-\lambda,\lambda]}(x).
\end{equation*}
For $\e<1$ we have
\begin{align*}
    & E[\sigma_\lambda]= \frac{\e}{4\lambda} - \frac{1}{8\lambda^2} \int_{-\lambda}^{\lambda} \int_{-\lambda}^{\lambda} G(x-y) dy dx\\
    & \ = \frac{\e}{4\lambda} - \frac{1}{4\lambda}\int_{-\lambda}^{\lambda} G(z) dz = \frac{1}{4\lambda}\left(\e - \int_{-\lambda}^{\lambda} G(z) dz\right)
\end{align*}
and since
\begin{equation*}
    \int_{-\lambda}^{\lambda} G(z) dz \rightarrow 1\quad \hbox{as}\ \ \lambda\rightarrow +\infty,
\end{equation*}
we easily obtain that there exists a $\overline{\lambda}$ such that $E[\sigma_{\overline{\lambda}}]<0$.
\endproof

\section{Stationary solutions in the $1$-$d$ case}\label{sec:stationary_one_d}

In this section we prove the main result of our paper, namely that nontrivial stationary solutions (which always exist in the case $\e<1$) in one space dimension with fixed mass and center of mass are unique. We shall first provide certain necessary conditions on the steady states and then prove that they are unique under such conditions. The main tool in this procedure is the use of the strong version of Krein-Rutman Theorem \ref{thm:KR2}.

We start with a necessary condition on the steady states which deals with a property of their support.

\begin{lemma}[Steady states have connected support]\label{lem:connected} Let $\rho$ be a stationary solution to \eqref{eq:main_stat_multid} in one space dimension, namely
\begin{equation}\label{eq:steady_oned}
    \rho\dep_x \left(\e \rho - G*\rho\right) =0\qquad \hbox{a.e. on}\ \R.
\end{equation}
Then, $\mathrm{supp}(\rho)$ is a connected set.
\end{lemma}

\proof
Let $\rho$ solve \eqref{eq:steady_oned}. Let us first assume that $\rho$ is compactly supported. Suppose that $\mathrm{supp}(\rho)$ is not connected. Accordingly, let $[a,b]$ be a non trivial interval such that
\begin{align}
    & \rho(x)\not\equiv 0& &\hbox{if}\ x<a&\nonumber\\
    & \rho(x) =0& &\hbox{if}\ a\leq x\leq b&\label{eq:supptho_twocomp}\\
    & \rho(x)\not\equiv 0& &\hbox{if}\ x>b.&\nonumber
\end{align}
Let us introduce the velocity field
\begin{equation*}
    V(x):=\begin{cases}
    -1 & \hbox{if}\ \ x\in(-\infty,a)\\
    1 & \hbox{if}\ \ x\in(b,+\infty)
    \end{cases}
\end{equation*}
and let $V\in C^1(\R)$. Let $u(x,s)$ be a local solution to the Cauchy problem for the continuity equation
\begin{equation*}
    \begin{cases}
        \dep_s u + \dep_x (u V) = 0& \\
        u(x,0)=\rho(x). &
    \end{cases}
\end{equation*}
Let us compute the evolution of the energy $E$ along $u$ at the time $s=0$:
\begin{align*}
   & \frac{d}{ds} E[u(s)] |_{s=0}= \int u_s(\e u(x,s) - G*u(x,s)) dx|_{s=0} = \int \rho V \dep_x (\e \rho - G*\rho) dx = 0.
\end{align*}
Then, by definition of $V$ we have
\begin{align}
    & \e \int_\R \rho V\dep_x \rho =\frac{\e}{2}\int_{-\infty}^a \dep_x \rho^2 dx +\frac{\e}{2}\int_a^b V \dep_x \rho^2 dx - \frac{\e}{2}\int_b^{+\infty}\dep_x \rho^2 dx =0\label{eq:use_comp_supp}
\end{align}
because $\rho_x=0$ on $[a,b]$ and $\rho=0$ on $x=a,b$ and at $\pm\infty$.
Hence, we have
\begin{equation}\label{lemma_support1}
    0 = \int_{-\infty}^{+\infty} \rho V \dep_x G*\rho dx = -\int_{-\infty}^a \rho G'*\rho dx + \int_b^{+\infty} \rho G'*\rho dx.
\end{equation}
We compute
\begin{align*}
   & \int_{-\infty}^a \rho G'*\rho dx =  \int_{-\infty}^a  \int_{-\infty}^a \rho(x)G'(x-y)\rho(y) dy dx +  \int_{-\infty}^a  \int_b^{+\infty} \rho(x)G'(x-y)\rho(y) dy dx,
\end{align*}
the first term on the above right-hand side is zero since $G'$ is odd and the integration domain is symmetric in $x$ and $y$. Since $G'(z)\geq0$ as $z\leq0$, we have for the second term
\begin{equation*}
   \rho(x)G'(x-y) \rho(y) \geq 0 \quad\hbox{on} \ (x,y)\in (-\infty,a)\times (b,+\infty).
\end{equation*}
In a similar way one can prove that
\begin{equation*}
    \int_b^{+\infty}\rho G'*\rho dx = \int_b^{+\infty}dx\int_{-\infty}^a dy\rho(x)G'(x-y) \rho(y)
\end{equation*}
with the integrand $\rho(x)G'(x-y) \rho(y)\leq 0$ on the integration domain. Therefore, \eqref{lemma_support1} implies that
\begin{equation}
    \rho(x)\rho(y) \equiv 0\qquad \hbox{on}\ \ \{x<a\}\cap\{y>b\}.\label{eq:supprho_condition}
\end{equation}
We have thus proven that, whenever \eqref{eq:supptho_twocomp} holds, then \eqref{eq:supprho_condition} has to be necessarily satisfied. Let $A, B$ be two nonempty connected components of $\mathrm{supp}(\rho)$ and let $[\alpha,\beta]$ be the maximal interval such that
\begin{align*}
    & a<b,\quad \hbox{for all}\ a\in A,\ b\in B\\
    & \rho\equiv 0\quad \hbox{on}\ [\alpha,\beta]\\
    & \alpha\geq a,\quad \hbox{for all}\ a\in A \\
     & \beta\leq b,\quad \hbox{for all}\ b\in B.
\end{align*}
Then, $\rho(x)\rho(y)=0$ for all $(x,y)$ such that
\begin{equation*}
    x<\alpha,\quad y>\beta,
\end{equation*}
which implies that either $A$ or $B$ cannot be in the support of $\rho$, and that is a contradiction. In order to generalize the proof to a stationary solution $\rho$ which is not compactly supported, one can cutoff $\rho$ to have compact support in such a way that the $L^2$ norm of the compactly supported approximation is arbitrarily close to the $L^2$ norm of $\rho$. Then, the estimate
\begin{equation*}
    \e\int|\rho\dep_x \rho| dx \leq \|G\|_{L^1}\|\rho\|_{L^2}^2
\end{equation*}
implies that the integrals in \eqref{eq:use_comp_supp} converge at infinity, therefore all the above computations are valid up to an arbitrary difference which vanishes in the limit.
\endproof

\begin{remark}
\emph{In the case $\mathrm{supp}(G)=[-R,R]$ one can use the same strategy as in Lemma \ref{lem:connected} to prove that, given two connected components $A,B$ of $\mathrm{supp}(\rho)$ one has $\mathrm{dist}(A,B)>2R$. The proof is a straightforward generalization of the above arguments, and it is therefore left to the reader.}
\end{remark}

We now exploit a standard symmetric rearrangement technique to prove that the minimizers of the energy are symmetric and monotonically decreasing on $x>0$ under the constraint of zero center of mass, cf. \cite{lieb,lieb_loss}.

\begin{prop}\label{prop:rearrange}
Let $\rho_\infty$ be a minimizer for the energy
\begin{equation*}
    E[\rho]=\frac{\e}{2}\int\rho^2 (x) dx - \frac{1}{2}\int\int G(x-y) \rho(x) \rho(y) dy dx
\end{equation*}
under the constraint that the center of mass is zero. Then, $\rho_\infty$ is symmetric and monotonically decreasing on $x>0$.
\end{prop}

\proof
We have to prove that the energy decreases strictly when a function $u$ which is not symmetric and decreasing on $x>0$ is replaced by its symmetric rearrangement
\begin{equation}
u^*(x)=\sup\left\{t\geq 0:\mathrm{meas}(\left\{u>t\right\})> 2\left|x\right|\right\}.
\end{equation}
For every exponent $p\geq 1$ the following holds:
\begin{equation}
\int_{\mathbb R}(u^*)^pdx=\int_{\mathbb R}(u)^pdx,
\label{eq:equim}
\end{equation}
therefore the $L^2$ part of the energy is invariant when passing from $u$ to $u^*$. As for the interaction energy, we recall the well known \emph{Riesz's rearrangement inequality}, see e. g. \cite{lieb_loss},
\begin{equation}\label{eq:riesz}
    \int_{\R^d}\int_{\R^d} f(x) g(x-y) h(y) dydx \leq \int_{\R^d} \int_{\R^d} f^*(x) g^*(x-y) h^*(y) dydx,
\end{equation}
which holds for all nonnegative functions $f,g,h$ vanishing at infinity. Moreover, if $g$ is strictly decreasing on $x>0$ and symmetric, then equality in \eqref{eq:riesz} holds if and only if $f(x)=f^*(x-x_0)$ and $h(x)=h^*(x-x_0)$ for some common $x_0$. Apply such a theorem to our case, using $G^*=G$ and the fact that $u$ is not symmetric up to translations.
\endproof

Let us rephrase Lemma \ref{lem:trick} in the one-dimensional case.
\begin{lemma}\label{lem:constant_is_energy}
Let $\rho$ be a $1d$ steady state, i. e.
\begin{equation*}
    \e\rho = G*\rho + C\qquad \hbox{on}\ \mathrm{supp}[\rho]
\end{equation*}
for some $C\in \R$. Then, $C=2E[\rho]$.
\end{lemma}

\proof
The support of $\rho$ is connected in view of Lemma \ref{lem:connected}, therefore Lemma \ref{lem:trick} applies.
\endproof

\begin{lemma}\label{lem:translation}
Let $\rho\in \pp\cap L^2$ and let $x_0\in\R$. Let $\rho_{x_0}$ be defined by
\begin{equation*}
    \rho_{x_0}(x):=\rho(x+x_0).
\end{equation*}
Then, $E[\rho_{x_0}]=E[\rho]$.
\end{lemma}

\proof
It follows by direct computation of the energy and by change of variable under the integral sign.
\endproof

\begin{lemma}\label{lem:support_bounded}
Let $\rho$ be a steady state with $\e<1$. Then, the support of $\rho$ is compact.
\end{lemma}

\proof
We know from Lemma \ref{lem:connected} that the support of $\rho$ is a connected set. Suppose that $\mathrm{supp}(\rho)$ is not bounded. That means that $\mathrm{supp}(\rho)$ is of the form $(-\infty,b)$ ($b$ possibly $+\infty$) or $(a,+\infty)$ ($a$ possibly ($-\infty$). Assume first $\mathrm{supp}(\rho)=(a,+\infty)$. Then, Lemma \ref{lem:constant_is_energy} implies
\begin{equation*}
    2E[\rho] = \e\rho(x) - \int_a^{+\infty} G(x-y)\rho(y) dy = 0
\end{equation*}
for all $x\in (a,+\infty)$. Now, there are two possibilities: either $a=-\infty$ or $a>-\infty$. In the latter case, evaluation on $x=a$ implies
\begin{equation*}
    0 = \e \rho(a) = \int_a^{+\infty} G(a-y)\rho(y) dy
\end{equation*}
which is a contradiction because the integral on the right hand side is strictly positive in view of $\mathrm{supp}(G)=\R$. In the former case $a=-\infty$ we have then $\mathrm{supp}(\rho)=\R$, which implies
\begin{equation*}
    \e \rho(x) = \int_{-\infty}^{+\infty}G(x-y) \rho(y) dy
\end{equation*}
for all $x\in \R$. We can therefore integrate over $\R$ to obtain
\begin{equation*}
    \e = \|G\|_{L^1}= 1
\end{equation*}
which is a contradiction. The same proof can be produced in the case $\mathrm{supp}(\rho)=(-\infty,b)$. Therefore, the support of $\rho$ can only be a bounded interval.
\endproof

\begin{lemma}\label{lem:steady_sym}
Let $\rho$ be a steady state. Then there exists a symmetric steady state $\rhot$ such that
\begin{equation*}
    E[\rhot] = E[\rho].
\end{equation*}
\end{lemma}

\proof
From Lemma \ref{lem:connected} and Lemma \ref{lem:support_bounded} we know that $\mathrm{supp}[\rho]=(a,b)$ for some $a,b\in\R$. For a given $x\in(a,b)$ we have
\begin{equation}\label{eq:stat_rho}
    \e\rho(x) = G*\rho ( x) + C
\end{equation}
for some $C\in \R$. Evaluation on $x=a$ and $x=b$ gives
\begin{equation*}
    C=-\int_a^b G(a-y)\rho(y) dy = -\int_a^b G(b-y)\rho(y) dy.
\end{equation*}
Let $\rhob(x)=\rho(x+x_0)$ with $x_0=(a+b)/2$. Then $\rhob$ is still a steady state and it satisfies $E[\rhob]=E[\rho]$ thanks to Lemma \ref{lem:translation}. Moreover, the support of $\rhob$ is symmetric. Let us introduce
\begin{equation*}
    \rhot(x):=\frac{1}{2}(\rhob(x) + \rhob(-x)).
\end{equation*}
Clearly, $\mathrm{supp}[\rhot]=\mathrm{supp}[\rhob]$ and we have, for all $x\in\mathrm{supp}[\rhot]$,
\begin{align*}
    & \e\rhot(x) =\frac{\e}{2}(\rhob(x) + \rhob(-x)) = \frac{\e}{2}(\rho(x+x_0) + \rho(-x + x_0))\\
    & \ = \frac{1}{2}\int_{a}^{b}G(x+x_0-y)\rho(y) dy + \frac{1}{2}\int_{a}^{b}G(-x+x_0-y)\rho(y) dy + C \\
    & \ = \frac{1}{2}\int_{(a-b)/2}^{(b-a)/2}G(x-z)\rhob(z) dy + \frac{1}{2}\int_{(a-b)/2}^{(b-a)/2}G(-x-z)\rhob(z) dy + C\\
    & \ = \frac{1}{2}\int_{(a-b)/2}^{(b-a)/2}G(x-z)\rhob(z) dy + \frac{1}{2}\int_{(a-b)/2}^{(b-a)/2}G(x-z)\rhob(-z) dy + C\\
    & \ = \int_{(a-b)/2}^{(b-a)/2}G(x-z)\frac{1}{2}\left(\rhob(z) + \rhob(-z)\right) dz + C = \int_{(a-b)/2}^{(b-a)/2}G(x-z)\rhot(z) dz + C
\end{align*}
where we have used the symmetry of $G$. The above computation shows that $\rhot$ has the same energy as $\rho$ in view of the results in Lemma \ref{lem:translation} and Lemma \ref{lem:constant_is_energy}.
\endproof

\begin{lemma}[Support of a minimizer]\label{lem:support_minimizer}
Let $\rho_\infty$ be a global minimizer to $E$. Let $\rho$ be a steady state such that
\begin{equation*}
    \mathrm{meas}(\mathrm{supp}[\rho_\infty])\leq\mathrm{meas}(\mathrm{supp}[\rho]).
\end{equation*}
Then $\rho$ is also a minimizer.
\end{lemma}

\proof
Assume first that we are in the special case $\mathrm{supp}[\rho_\infty]\subseteq\mathrm{supp}[\rho]$. Let us compute the second variation of $E$ around the minimizer $\rho_\infty$ along the direction $\rho_\infty - \rho$.
\begin{align*}
    & \frac{d^2}{d\delta^2} E[\rho_\infty + \delta (\rho - \rho_\infty)]|_{\delta =0} \\
    & \ = \e\int (\rho - \rho_\infty)^2 dx - \int\int G(x-y)(\rho(x) - \rho_\infty(x))(\rho(y) - \rho_\infty(y))dy dx\\
    & \ = 2 E[\rho] + 2 E[\rho_\infty] - 2\int_{\mathrm{supp}[\rho_\infty]}\rho_\infty(\e\rho - G*\rho) dx\\
    & \ = 2 E[\rho] + 2 E[\rho_\infty] - 4 E[\rho]
\end{align*}
where the last step is justified by the fact that $\mathrm{supp}[\rho]\subseteq\mathrm{supp}[\rho_\infty]$. Therefore, since $\rho_\infty$ is a minimizer, the second derivative above is nonnegative, i. e.
\begin{equation*}
   0\leq  \frac{d^2}{d\delta^2} E[\rho_\infty + \delta (\rho - \rho_\infty)]|_{\delta =0} = 2 (E[\rho_\infty] - E[\rho]),
\end{equation*}
which yields $E[\rho]\leq E[\rho_\infty]$. Since $\rho_\infty$ is a minimizer, then so is $\rho$.
In the general case in which $\mathrm{supp}[\rho_\infty]\nsubseteq\mathrm{supp}[\rho]$, consider a translation $\rho_{x_0}(x)=\rho(x-x_0)$ in such a way that the support of $\rho$ contains the support of $\rho_\infty$. Since the energy is invariant after translation in view of Lemma \ref{lem:translation}, the assertion is proven.
\endproof

%


We are now getting closer to the proof of our uniqueness result. Let us recall the following important theorems, see e. g. \cite{KR} and the references therein.

\begin{thm}
[Krein--Rutman Theorem]\label{thm:KR}
Let $X$ be a Banach space, let $K\subset X$ be a \emph{total cone}, i. e. such that $\lambda K\subset K$ for all $\lambda \geq 0$ and such that the set $\{u-v,\ u,v\in K\}$ is dense in $X$. Let $T$ be a compact linear operator such that $T(K)\subset K$ with positive spectral radius $r(T)$. Then $r(T)$ is an eigenvalue for $T$ with an eigenvector $u\in K\setminus \{0\}$.
\end{thm}

An important consequence \cite{KR} of the Krein--Rutman theorem, which will be extremely useful in the sequel, is the following
\begin{thm}[Krein--Rutman Theorem, strong version]\label{thm:KR2}
Let $X$ be a Banach space, $K\subset X$ a \emph{solid cone}, i. e. such that $\lambda K\subset K$ for all $\lambda \geq 0$ and such that $K$ has a nonempty interior $K_0$. Let $T$ be a compact linear operator which is \emph{strongly positive} with respect to $K$, i. e. such that $T[u]\in K_0$ if $u\in K$. Then,
\begin{enumerate}
  \item [(i)] The spectral radius $r(T)$ is strictly positive and $r(T)$ is a simple eigenvalue with an eigenvector $v\in K_0$. There is no other eigenvalue with a corresponding eigenvector $v\in K$.
  \item [(ii)] $|\lambda|<r(T)$ for all other eigenvalues $\lambda\neq r(T)$.
\end{enumerate}
\end{thm}

We shall now prove the uniqueness of symmetric steady states with unit mass which are monotonically decreasing on the positive semi-axis in the case $\varepsilon<1$. We already know that under the above assumptions we can write, for $x\in[-L,L]=\mathrm{supp}[\rho]$,
\begin{equation}\label{eq:stat_compact}
    \varepsilon \rho (x) = \int_{-L}^{L} G(x-y) \rho(y) dy + C,\qquad C=2 E[\rho].
\end{equation}
Taking the derivative w.r.t $x\in [-L,L]$ we obtain
\begin{equation*}
    \varepsilon \rho' (x) = \frac{d}{dx}\int_{-L}^{L} G(x-y) \rho(y) dy = \frac{d}{dx} G*\rho (x) = \int_{-L}^{L} G(x-y) \rho'(y) dy.
\end{equation*}
The symmetry of $\rho$ and $G$ implies, for $x\in[0,L]$,
\begin{equation*}
     \varepsilon \rho' (x) = - \int_{0}^{L} G(x+y) \rho'(y) dy + \int_{0}^{L} G(x-y) \rho'(y) dy = \int_0^L \left[G(x-y)- G(x+y)\right] \rho'(y) dy.
\end{equation*}
Assuming that $\rho\in C^1([-L,L])$, finding a steady state with the above assumptions is equivalent to find $\rho$ on $[0,L]$ such that
\begin{align*}
    & \rho(L)=0,\\
    & -\rho'(x)=u(x),\quad x\in[0,L],\\
    & u\geq 0,\quad \hbox{and $u$ solves}\quad \varepsilon u = \int_0^L H(x,y) u(y) dy,\\
    & H(x,y)=G(x-y)- G(x+y).
\end{align*}
To convince ourselves about that, integrate
\begin{equation*}
    -\e\rho'(x) = -\int_0^L (G(x-y) - G(x+y))\rho'(y) dy
\end{equation*}
over the interval $[\xi,L]$ for some $\xi\in[0,L)$. Then $\rho(L)=0$ and integration by parts imply
\begin{align*}
     -\e\rho(\xi) = & - \int_\xi^L dx \int_0^L (G(x-y) - G(x+y))\rho'(y) dy \\
     = & -\int_\xi^L dx [(G(x-L) - G(x+L))\rho(L)-(G(x) - G(x))\rho(0)]\\
      & + \int_\xi^L dx \int_0^L (-G'(x-y) - G'(x+y))\rho(y) dy \\
    = & \int_0^L \rho(y) dy \int_\xi^L (-G'(x-y) - G'(x+y)) dx\\
     = & \int_0^L \rho(y) [ -G(L-y) - G(L+y) + G(\xi-y) + G(\xi+y) ]
\end{align*}
which implies, by the symmetry of $G$,
\begin{equation*}
    \e\rho(x) = \int_{-L}^L G(x-y) \rho(y) dy + C,\qquad C = -\int_{0}^{L}(G(L-y) + G(L+y)) \rho(y) dy.
\end{equation*}
For further reference, we introduce the operator
\begin{equation}\label{eg:Gl}
    \mathcal{G}_L[\rho](x) :=\int_0^L\left[G(x-y) + G(x+y) - G(L-y) - G(L+y)\right]\rho(y) dy
\end{equation}
on the Banach space
\begin{equation*}
    Y_L:=\left\{\rho\in C([0,L])\ :\ \rho(L)=0\right\}.
\end{equation*}
In order to simplify the notation, we also define the following operator
\begin{equation*}
    \mathcal{H}_L[u] (x) :=\int_0^L H(x,y) u(y) dy = \int_0^L (G(x-y) - G(x+y))u(y) dy.
\end{equation*}

\begin{prop}\label{prop:KRderivatives}
For a fixed $L>0$ there exists a unique symmetric function $\rho\in C^2([-L,L])$ with unit mass and with $\rho'(x)\leq 0$ on $x\geq 0$ such that $\rho$ solves
\eqref{eq:stat_compact} for some $\e=\varepsilon(L)>0$. Such function $\rho$ also satisfies $\rho''(0)<0$. Moreover, $\e(L)$ is the largest eigenvalue of the compact operator $\mathcal{G}_L$ on the space Banach $Y_L$ and any other eigenfunction of $\mathcal{G}_L$ on $Y_L$ with unit mass has the corresponding eigenvalue $\e'$ satisfying $|\e'|<\e(L)$.
\end{prop}

\proof
Since $G$ is decreasing on the half-line $[0,+\infty)$ we get
\begin{equation*}
    H(x,y)= G(x-y)- G(x+y)\geq 0,\qquad \hbox{on}\quad x,y \geq 0.
\end{equation*}
Consider now the Banach space
\begin{equation*}
    X_L=\left\{ f\in C^1([0,L]):\ f(0)=0\right\}
\end{equation*}
endowed with the $C^1$ norm
\begin{equation*}
    \|f\|_{X_L} =\|f\|_{L^\infty([0,L])} + \|f'\|_{L^\infty([0,L])}.
\end{equation*}
It can be easily seen that the set
\begin{equation*}
    K:=\left\{f\in X:\ f \geq 0\right\}
\end{equation*}
is a solid cone in $X$. Indeed, any function $f\in K$ with $f'(0)>0$ is in the interior of $K$. Moreover, for a given $u\in K$, we have
\begin{equation*}
   \mathcal{H}_L[u](x)= \int_0^L H(x,y) u(y) dy \geq 0
\end{equation*}
for all $x\in[0,L]$ and
\begin{equation*}
    \mathcal{H}_L[u](0) = \int_0^L H(0,y) u(y) dy = \int_0^L (G(-y) - G(y))u(y) dy =0.
\end{equation*}
Therefore $\mathcal{H}$ is a positive operator in the sense provided by the definition of the cone $K$. Indeed, we can prove that $\mathcal{H}$ is strongly positive, i. e. for a given $u\in K$, $\mathcal{H}[u]$ belongs to the interior of $K$. In order to see that, for a $u\in K\setminus\{0\}$ compute
\begin{equation*}
    (\mathcal{H}_L[u])'(0)= \int_0^L ( G'(-y) - G' (y)) u(y) dy =-2\int_0^L G'(y) u(y) dy >0,
\end{equation*}
and therefore $\mathcal{H}_L[u]$ belongs to the interior of $K$. Hence, we can apply the stronger version of Krein--Rutman theorem \ref{thm:KR2}, which implies the existence of a simple eigenvalue $\varepsilon>0$ equal to the spectral radius of $\mathcal{H}_L$. More precisely, there exists a family of solutions $u$ to \begin{equation*}
    \varepsilon u= \mathcal{H}_L[u]
\end{equation*}
generated by one given nontrivial element $\bar u$ in the interior of $K$. This implies that the corresponding set of symmetric and monotone $\rho$ solving \eqref{eq:stat_compact} satisfies
\begin{equation*}
    \rho (x) = \rho (x) - \rho (L) = - \int_x^L \rho'(y) dy = \int_x^L u(y) dy = \alpha \int_x^L \bar u (y) dy,
\end{equation*}
for $\alpha>0$. We choose $\alpha$ as follows
\begin{equation*}
    \alpha = \left(2\int_0^L \int _x^L \bar u(y) dy dx\right)^{-1},
\end{equation*}
and we obtain that $\rho$ has unit mass on $[-L,L]$. It is clear that $\rho'(x)\leq 0$ for $x\geq 0$, $\rho'(0)=0$, and $\rho''(0)<0$. In view of the statement (i) of theorem \ref{thm:KR2}, there exists no other eigenvalues to $\mathcal{H}_L$ with eigenvectors in $K$ besides the one $\e$ with eigenfunction $\bar u$, and all other eigenvalues $\e'$ with eigenfunctions in $X_L$ satisfy $|\e'|<\e$.
\endproof

The eigenvalue $\varepsilon$ (which coincides with the spectral radius of $\mathcal{H}_L$) can be considered as a function of $L$, namely $\varepsilon=\varepsilon(L)$. The behavior of such function is established in the next proposition.

\begin{prop}[Behavior of the function $\e(L)$]\label{prop:epsilonL}
The simple eigenvalue $\varepsilon(L)$ found in Proposition \ref{prop:KRderivatives} is uniquely determined as a function of $L$ with the following properties
\begin{itemize}
  \item [(i)] $\varepsilon(L)$ is strictly increasing with respect to $L$
  \item [(ii)] $\lim_{L\rightarrow +\infty}\varepsilon (L)= 1$
  \item [(iii)] $\varepsilon(0)=0$.
\end{itemize}
\end{prop}

\proof
In order to prove (i), let us consider the equation
\begin{equation*}
    \varepsilon(L) u_L(x) = \mathcal{H}_L[u_L](x) = \int_0^L H(x,y) u_L(y) dy,\qquad x\in[0,L],
\end{equation*}
where $u_L$ is the unique eigenfunction obtained in Proposition \ref{prop:KRderivatives}.
We multiply the above equation by $u_L(x)$ and integrate over $[0,L]$ to obtain
\begin{equation*}
    \varepsilon(L)\int_0^L u_L(x)^2 dx = \int_0^L \mathcal{H}_L[u_L](x) u(x)dx.
\end{equation*}
Recall that the eigenvalue $u_L$ satisfies $u_L(0)=0$ and, for $x\in(0,L]$,
\begin{equation*}
    u_L(x)=\frac{1}{\e(L)}\int_0^L H(x,y) u(y) dy >0
\end{equation*}
since $H(x,y)=G(x-y)-G(x+y)>0$ for all $y\in[0,L]$ under the assumption $x>0$ in view of the strict decreasing monotonicity of $G$ on $x>0$.
For a general $L\in (0,+\infty)$ and a $\delta>0$ (small enough) we have
\begin{align}
    & I_1:=\e(L+\delta)\int_0^{L+\delta} u_{L+\delta}^2(x) dx - \e(L)\int_0^{L} u_{L}^2(x) dx\nonumber\\
    & \ = \int_0^{L+\delta} \mathcal{H}_{L+\delta}[u_{L+\delta}](x) u_{L+\delta}(x)dx - \int_0^{L} \mathcal{H}_{L}[u_{L}](x) u_{L}(x)dx=:I_2.\label{epsilonincreasing1}
\end{align}
We analyze the two terms $I_1$ and $I_2$ separately. $I_1$ can be expanded as follows:
\begin{align*}
    & I_1 = (\e(L+\delta)-\e(L))\int_0^{L+\delta} u_{L+\delta}^2(x) dx\\
    & \  + \e(L)\int_0^{L+\delta}(u_{L+\delta}(x)-u_L(x))(u_{L+\delta}(x)+u_L(x))dx + \e(L)\int_L^{L+\delta} u_L^2(x) dx.
\end{align*}
$I_2$ is given by
\begin{align*}
    & I_2= \int_0^{L+\delta} \left(\mathcal{H}_{L+\delta}[u_{L+\delta}](x) - \mathcal{H}_{L}[u_{L}](x)\right)u_{L+delta}(x) dx \\
    & \ + \int_L^{L+\delta} \mathcal{H}_L[u_L](x)u_L(x) dx + \int_0^{L+\delta}\mathcal{H}_L[u_L](x)\left(u_{L+\delta}(x)-u_L(x)\right) dx\\
    & \ = \int_0^{L+\delta} \left(\mathcal{H}_{L+\delta}[u_{L+\delta}](x) - \mathcal{H}_{L}[u_{L}](x)\right)u_{L+delta}(x) dx \\
    & \ + \e(L)\int_L^{L+\delta}u_L^2(x) dx + \e(L)\int_0^{L+\delta}u_L(x)\left(u_{L+\delta}(x)-u_L(x)\right) dx
\end{align*}
and on substituting $I_1$ and $I_2$ in \eqref{epsilonincreasing1} we can cancel some terms and obtain
\begin{align*}
    & (\e(L+\delta)-\e(L))\int_0^{L+\delta} u_{L+\delta}^2(x) dx + \e(L)\int_0^{L+\delta} (u_{L+\delta}(x)-u_L(x)) u_{L+\delta}(x) dx\\
    & \ = \int_0^{L+\delta} \left(\mathcal{H}_{L+\delta}[u_{L+\delta}](x) - \mathcal{H}_{L}[u_{L}](x)\right)u_{L+\delta}(x) dx\\
    & \ = \int_0^{L+\delta}\int_L^{L+\delta} H(x,y)u_L(\delta)(y)u_{L+\delta}(x) dy dx \\
    & \ + \int_0^{L+\delta}\int_0^{L+\delta}H(x,y)(u_{L+\delta}(y)-u_L(y))u_{L+\delta}(x) dy dx\\
    & \ = \int_0^{L+\delta}\int_L^{L+\delta} H(x,y)u_L(y)u_{L+\delta}(x) dy dx \\
    & \ + \e(L+\delta)\int_0^{L+\delta}u_{L+\delta}(x)(u_{L+\delta}(x) - u_{L}(x))dx
\end{align*}
where we have used the definition of $\mathcal{H}_L$ and the property $H(x,y)=H(y,x)$. By suitably expanding the term on the left hand side in the above identity, we obtain
\begin{align*}
    & (\e(L+\delta) - \e(L))\int_0^{L+\delta} u_L(x) u_{L+\delta}(x) dx = \int_0^{L+\delta} \int_L^{L+\delta} H(x,y) u_L(y) u_{L+\delta}(x) dy dx\\
    & \ = \e(L+\delta)\int_L^{L+\delta}u_{L+\delta}(y)u_L(y) dy
\end{align*}
and the positivity property of $u_L$ implies that
\begin{equation*}
    \e(L+\delta)>\e(L),
\end{equation*}
which proves (i).


Let us now prove (ii). Assume by contradiction that
\begin{equation*}
    \lim_{L\rightarrow +\infty}\varepsilon (L)= \e_0<1.
\end{equation*}
Let $\e\in (\e_0,1)$. We know from Theorem \ref{thm:existence_minimizers} that there exists a minimizer $\rho_\e$ for the energy $E$ with zero center of mass. We also know that the support of $\rho_\e$ is compact from Lemma \ref{lem:support_bounded}. From Proposition \ref{prop:rearrange} we know that $\rho_\e$ is symmetric and monotonically decreasing on $x>0$. Therefore, $\rho_\e$ is the unique eigenfunction with unit mass provided by Proposition \ref{prop:KRderivatives}, and the support of $\rho_\e$ is $[-L,L]$ for some $L>0$. Therefore, the corresponding eigenvalue should be $\e(L)<\e_0$, which is a contradiction since $\e$ and $\e_0$ are two different eigenvalues with the same eigenfunction.

Let us prove (iii). By letting $L\searrow 0$ one has that the operator $\mathcal{H}_L$ is the zero operator, and therefore $\e(0)$ should be the eigenvalue of the zero operator, which can only be zero.

\endproof

We are now ready to prove the main result of this paper.

\begin{thm}\label{thm:main1d} Let $\e<1$. Then, there exists a unique $\rho\in L^2$ solution to
\begin{equation*}
    \rho\partial_x(\e \rho - G*\rho)=0,
\end{equation*}
with unit mass and zero center of mass. Moreover,
\begin{itemize}
  \item $\rho$ is symmetric and monotonically decreasing on $x>0$,
  \item $\rho\in C^2(\mathrm{supp}[\rho])$,
  \item $\mathrm{supp}[\rho]$ is a bounded interval in $\R$,
  \item $\rho$ has a global maximum at $x=0$ and $\rho''(0)<0$,
  \item $\rho$ is the global minimizer of the energy $E[\rho]=\frac{\e}{2}\int \rho^2 dx - \frac{1}{2}\int \rho G*\rho dx$.
\end{itemize}
\end{thm}

\proof
We know from Theorem \ref{thm:existence_minimizers} that there exists a minimizer $\rho_\infty$ with unit mass and zero center of mass, which is symmetric and monotonically decreasing on $x>0$ in view of Proposition \ref{prop:rearrange} and compactly supported on a certain $[-L,L]$ in view of Lemma \ref{lem:support_bounded}. From the results in Propositions \ref{prop:KRderivatives} and \ref{prop:epsilonL}, we know that there exists a unique steady state with such properties, because the correspondence $\e=\e(L)$ is one-to-one. So, the only possibility to violate uniqueness of steady states with unit mass and zero center of mass is to have a steady state which violates either the monotonicity property or the symmetry.
Suppose first that there exists a steady state with zero center of mass $\rho$ which is not symmetric, it is not restrictive to assume the support of $\rho$ is $[-L',L']$. Then, we know from Lemma \ref{lem:steady_sym} that it is possible to construct a symmetric steady state $\rhot$ with the same energy of $\rho$ and with the same support of $\rho$. Now, there are two possibilities: either $\rhot$ is a minimizer or not. In the former case $\rho$ is also a minimizer and this is a contradiction (a minimizer is symmetric). In the latter case, the support of $\rhot$ is strictly contained in the support of $\rho_\infty$ in view of Lemma \ref{lem:support_minimizer}, and $\rhot$ is not monotonically decreasing on $x>0$ because otherwise it would be the unique minimizer provided before. Therefore, with the notation of Proposition \ref{prop:KRderivatives}, $-\rhot'$ is an eigenfunction for $\mathcal{H}_{L'}$ in the space $X_{L'}$ which is not belonging to the solid cone $K$. Therefore, the stronger version of Krein-Rutman's Theorem \ref{thm:KR2} and the fact that $\e(L)$ is increasing imply that $L'>L$, since $\rho$ is an eigenfunction outside the solid cone $K$, and it therefore should have an eigenvalue strictly less than $\e(L')$. This implies that $\e(L)<\e(L')$ and therefore $L<L'$. Now this is a clearly a contradiction because we said before that the support or $\rhot$ is strictly contained in the support of $\rho_\infty$, so $L>L'$.
The case in which $\rho$ is symmetric but not monotonic on $x>0$ can be covered by repeating the same argument above (assume $\rho=\rhot$!).
\endproof

\begin{cor}[Concavity of $\rho$ for small $\e$] \label{cor:Concavity}
There exists a value $\e_0\in(0,1)$ such that, for all $\e\in(0,\e_0)$ the corresponding stationary solution provided in Theorem \ref{thm:main1d} is concave on the whole interval $[0,L]$.
\end{cor}

\proof
We can differentiate twice w.r.t $x$ in
\begin{equation*}
    \e\rho(x)=\int_{-L}^{L}G(x-y)\rho(y)dy + C
\end{equation*}
to obtain
\begin{equation*}
    \e\rho''(x)=\int_{-L}^{L}G''(x-y)\rho(y)dy
\end{equation*}
for all $x\in[-L,L]$. Therefore, $G''$ is evaluated on the interval $[-2L,2L]$ in the above integral. We know from Proposition \ref{prop:epsilonL} that $L$ is a monotonically increasing function of $\e$ with $\lim_{\e\searrow 0} L(\e)=0$. Since $G''(0)<0$, and $G\in C^2$, then there exists $L_0>0$ such that $G''<0$ on $[-2L_0,2L_0]$. Let $\e_0$ be the eigenvalue in $K$ corresponding to $L=L_0$. Then, the eigenfunction $\rho$ is concave on its support.
\endproof

\begin{remark}[The case $\mathrm{meas}(\mathrm{supp}(G))<+\infty$]
\emph{It is also interesting to consider the case with the support of $G$ being bounded, i.e. $\mathrm{supp}(G)=[-g,g]$ with $G$ being symmetric and monotone on $[-g,0]$. Then most techniques of the paper remain valid, in particular when considering simply connected stationary states, whose existence can be shown along the lines of the above arguments. With the same arguments as in the proof of Lemma \ref{lem:connected} one can show that the distance between two connected components of a stationary solution is at least $2g$. This means that stationary solutions consist of a countable number of connected components, which are not influenced by any other connected component (the distance is larger than the kernel range). Thus, each connected component is a stationary solution by itself. Thus, in particular we can prove the existence of  infinitely many stationary states.
The remaining step would be to characterize uniquely the behavior of each connected component as an energy minimizer. All arguments above apply to the case of simply connected solutions in the case of a finite range kernel $G$, except the application of the strong version of the Krein-Rutman theorem. If we cannot guarantee that the support of the global minimizer is contained in $[-\frac{g}2,\frac{g}2]$, there can be elements in the nonnegative cone such that the convolution with $G$ is not strictly positive. Although we strongly believe that such a smallness of the support holds at least for $\epsilon$ small, we so far did not succeed in proving such a result.}
\end{remark}

\section{Numerical Results}\label{sec:numerics}
In the following we show numerical simulations for the evolution equation \eqref{eq:main_evo}. We discretize the equation using an explicit Euler scheme and finite difference methods. In one dimension, we partitionate the domain $\Omega=[a,b]$ using an equidistant grid with $n+1$ grid points $a=x_0<x_1<...<x_n=b$ and step size $h=(b-a)/(n+1)$. Furthermore we use the following finite difference scheme:
\[\frac{\rho_i^{j+1}-\rho_i^{j}}{dt}=D_-^x(\rho_i^j D_+^x(\epsilon\rho_i^j-G\ast\rho_i^j))\]
with forward and backward difference quotients
\[D_+^x\rho_i=\frac{\rho_{i+1}- \rho_{i}}{h},\quad D_-^x\rho_i=\frac{\rho_{i}- \rho_{i-1}}{h}.\]
The time step size $dt$ has to be chosen appropriately, to guarantee stability. 
In a first example we consider an interaction potential $G(x) = \frac {1}{\sigma\sqrt{2\pi}}\exp\left(-\frac {1}{2} \left(\frac{x-\mu}{\sigma}\right)^2\right)$ with mean $\mu=0$ and variance $\sigma=1$, which fulfill the conditions (1)-(8),  $\left\|G\right\|_{ L^1}=1$  and models a wide range attraction.
\begin{figure}
			\subfigure[]{	\includegraphics[width=0.4\textwidth]{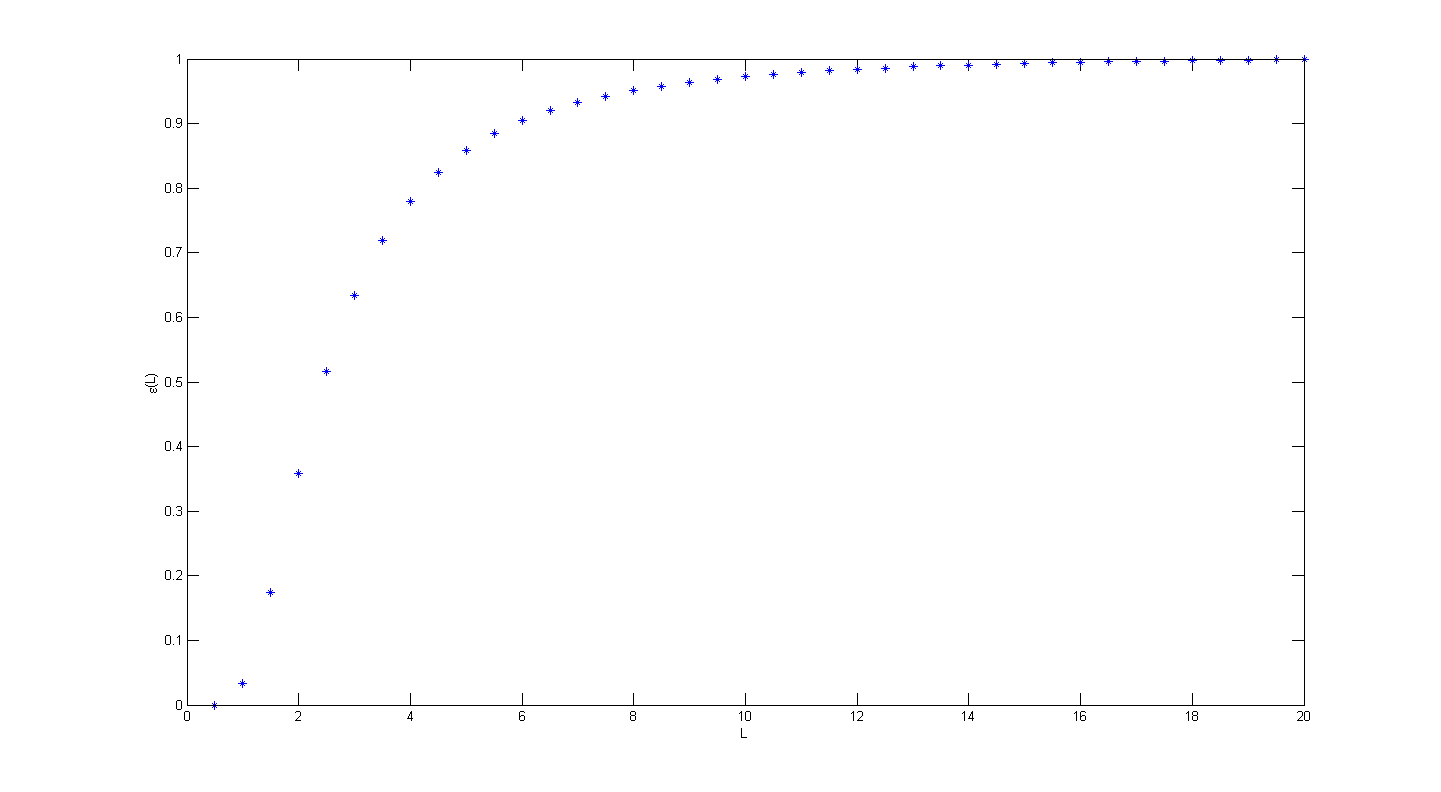}} \\
			\subfigure[]{ \includegraphics[width=0.3\textwidth]{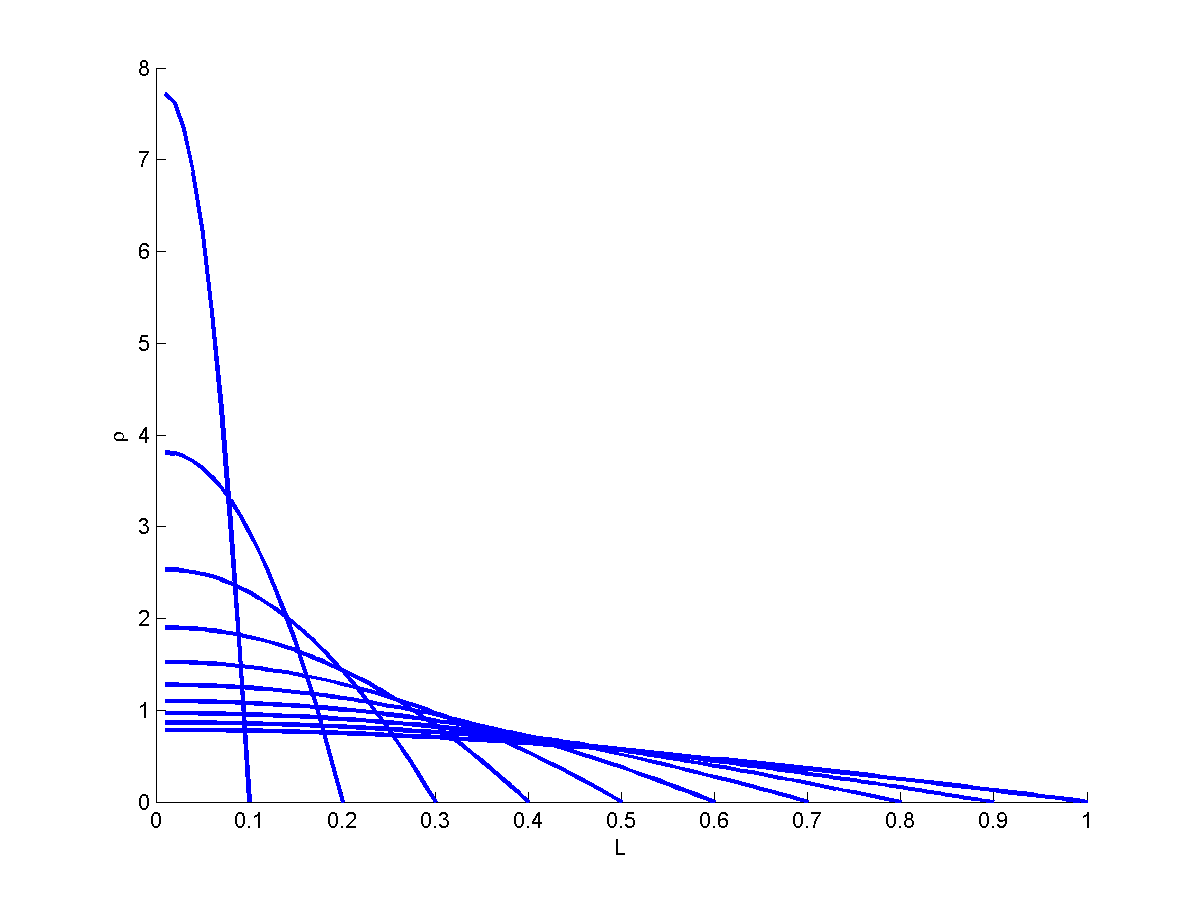}}
	\subfigure[]{	\includegraphics[width=0.3\textwidth]{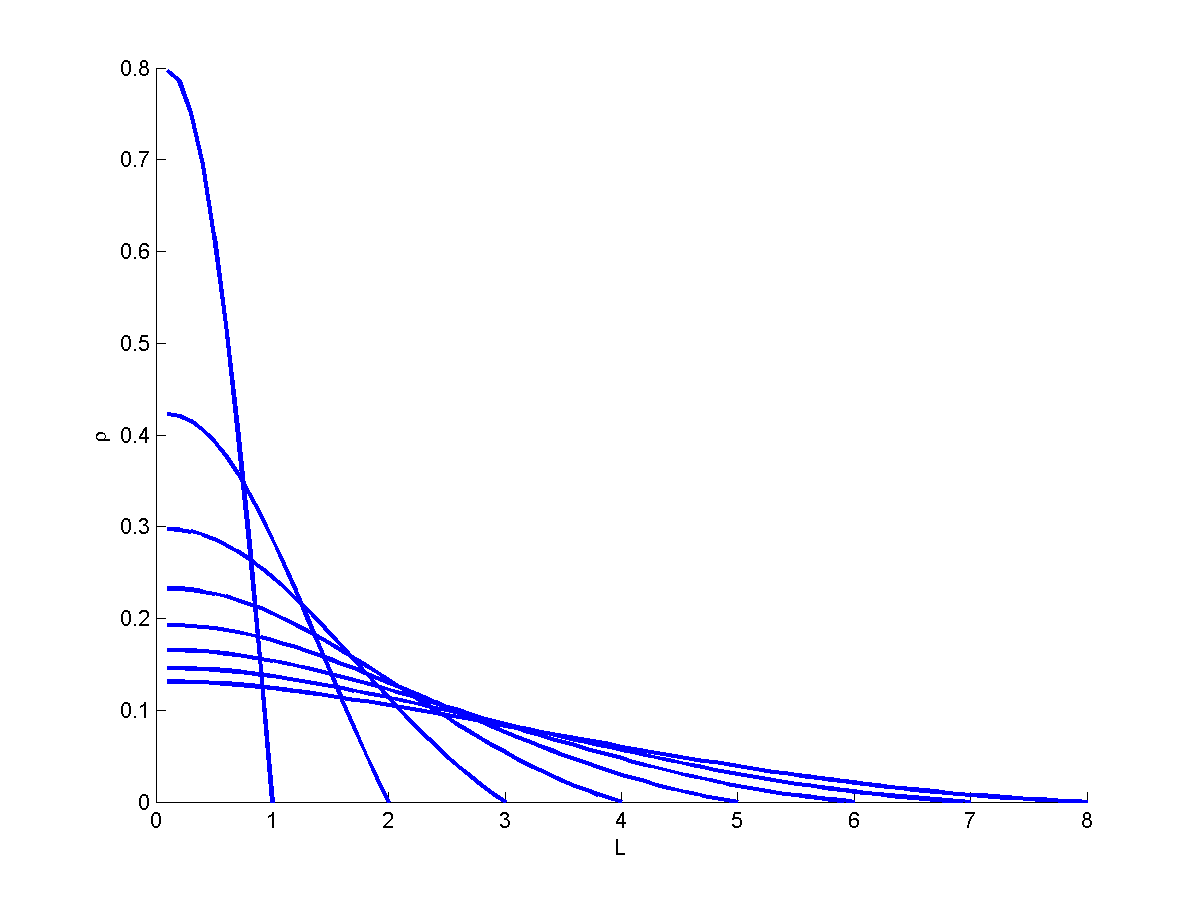}}
	\caption{Largest eigenvalues and corresponding eigenfunctions of the operator  $\mathcal{G}_L$ \ref{eg:Gl}: (a) Largest eigenvalues $\epsilon=\epsilon(L)$ of $\mathcal{G}_L$ on $L=(0,20]$; (b) Corresponding eigenfunctions for $\epsilon(L)$  with $L\in [0,1]$; (c) Corresponding eigenfunctions for $\epsilon(L)$  with  $L\in [1,8]$.}
		\label{fig:eigenfunctions}
\end{figure}
For this kernel we present in Figure \ref{fig:eigenfunctions} the solutions for the stationary equation \eqref{eq:main_stat}, which means we  calculated the largest eigenvalues  and  corresponding eigenvectors of the operator $\mathcal G_L$ defined in \eqref{eg:Gl} for different $L$.
The largest eigenvalues $\epsilon=\epsilon(L)$ are presented in Figure \ref{fig:eigenfunctions} (a). Like mentioned in Proposition \ref{prop:epsilonL}  $\epsilon(L)$ is strictly increasing with respect to $L$  and furthermore $\lim_{L\rightarrow \infty}\epsilon(L)=1$. The corresponding eigenfunctions with unit mass are presented in Figure \ref{fig:eigenfunctions} (b) and (c). We proved in Corollary \ref{cor:Concavity} the concavity of $\rho$ for small $\epsilon$. To better clarify the situation, we illustrate in Figure \ref{fig:eigenfunctions} (b) the eigenfunctions for $L\in (0,1]$, that means for $\epsilon <0.05$.  For a certain $\epsilon$,  which depends on the concavity of the kernel $G$,  the solution is not fully concave on its support any more, but bell shaped, compare Figure \ref{fig:eigenfunctions} (c).

\begin{figure}
	\centering
		\includegraphics[width=.45\textwidth]{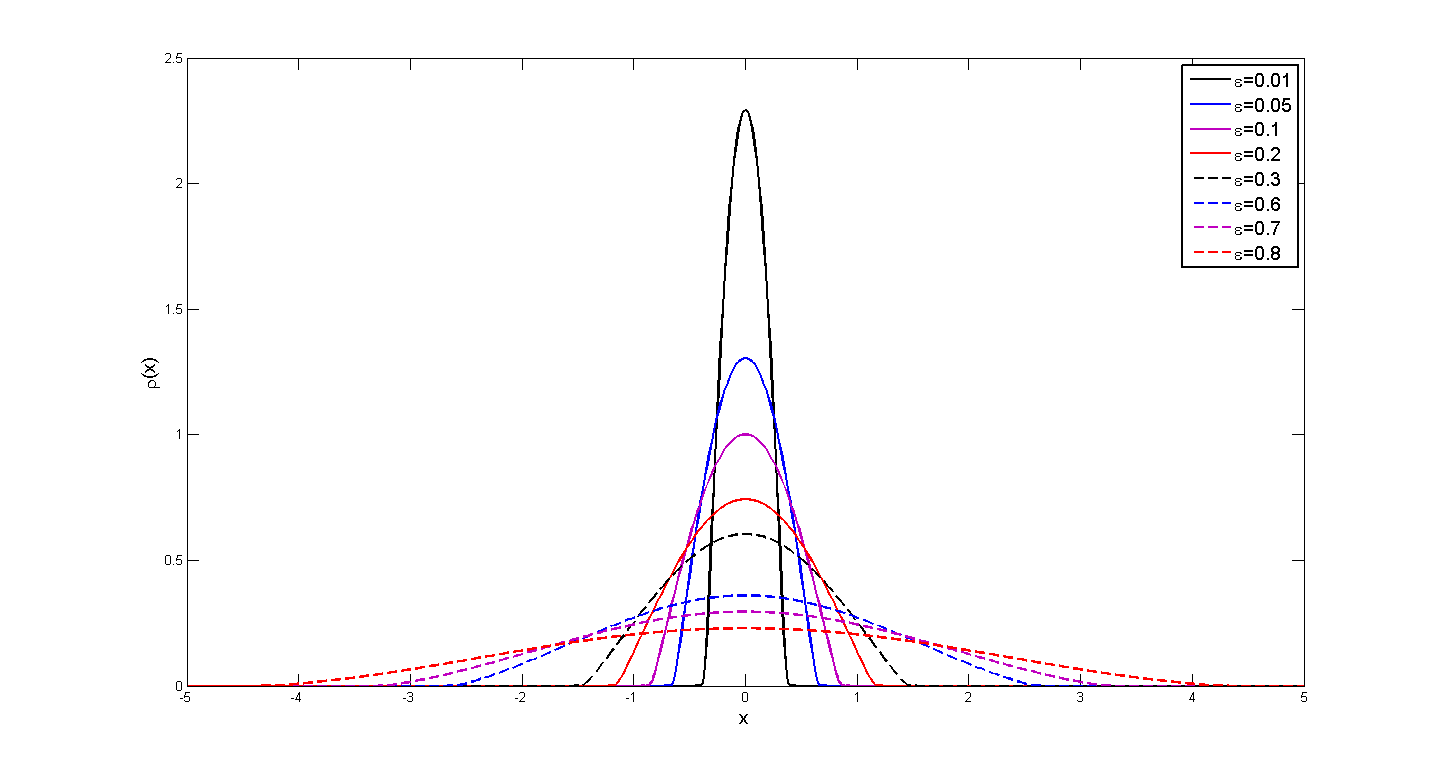}
\caption{Stationary solutions for equation \eqref{eq:main_evo} with $\epsilon \in (0,1)$.}
\label{fig:stationarySolutions1}
\end{figure}

To make this result more clear we present in Figure \ref{fig:stationarySolutions1} the stationary solutions of the evolution equation \eqref{eq:main_evo} for $\epsilon\in (0,1)$. We consider a compactly supported initial datum $\rho(x,0)=\rho_0$ with unit mass  $\int_{\Omega}\rho_0=1$. As the results are the same as in Figure \ref{fig:eigenfunctions} (b)-(c), we recognize again that up to a certain $\epsilon$ the solutions are concave and then  bell shaped. Furthermore we have mass conservation.
As proven before in this paper, for $\epsilon\geq \int G=1$ we do not have steady states because the impact of the diffusive term is higher then of the aggregation term. In this case we expect the solutions to behave like the self-similar Barenblatt-Pattle profiles. For $\epsilon=0$ we obtain an unique stationary solution (with zero center of mass), which is a Dirac-$\delta$-distribution with unit mass centered at zero.

\begin{figure}
	\subfigure[]{ \includegraphics[width=0.45\textwidth]{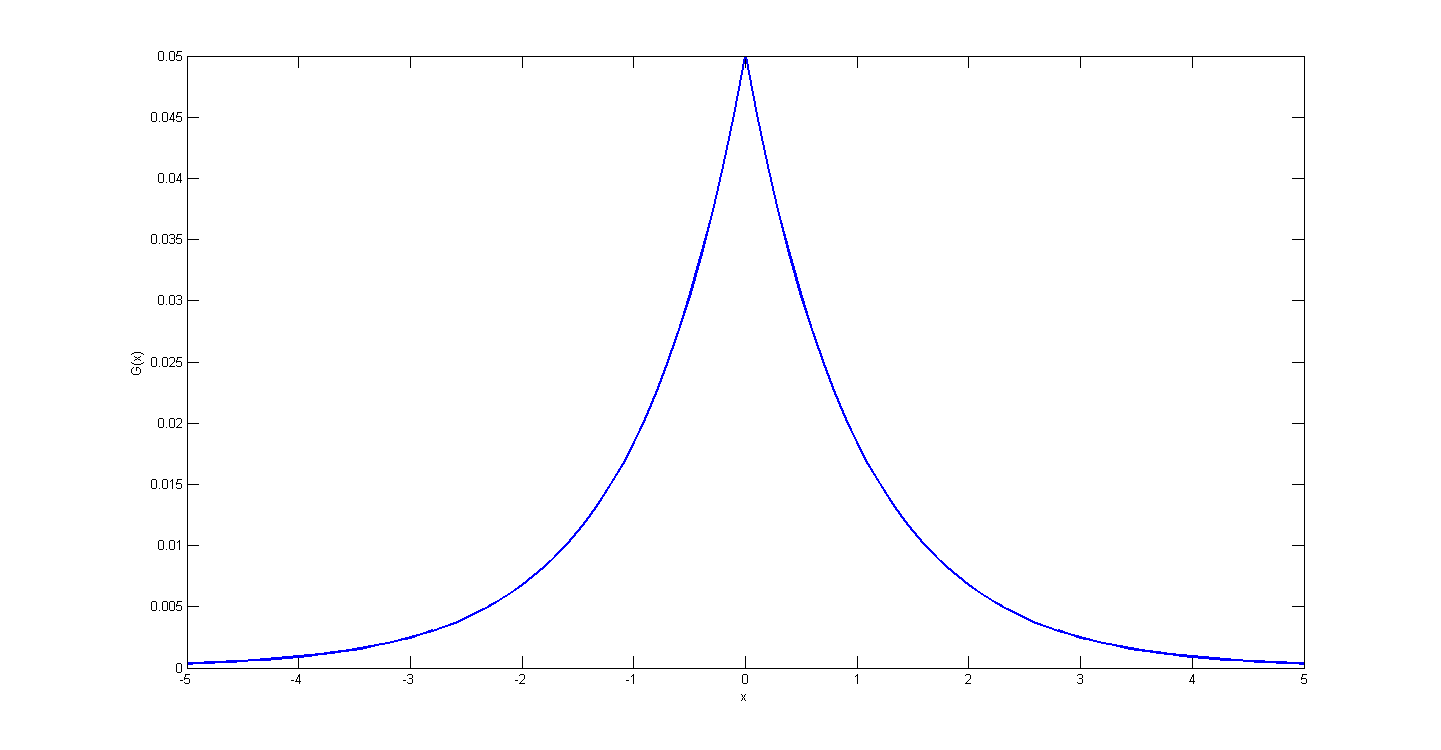}}
				\subfigure[]{ \includegraphics[width=0.45\textwidth]{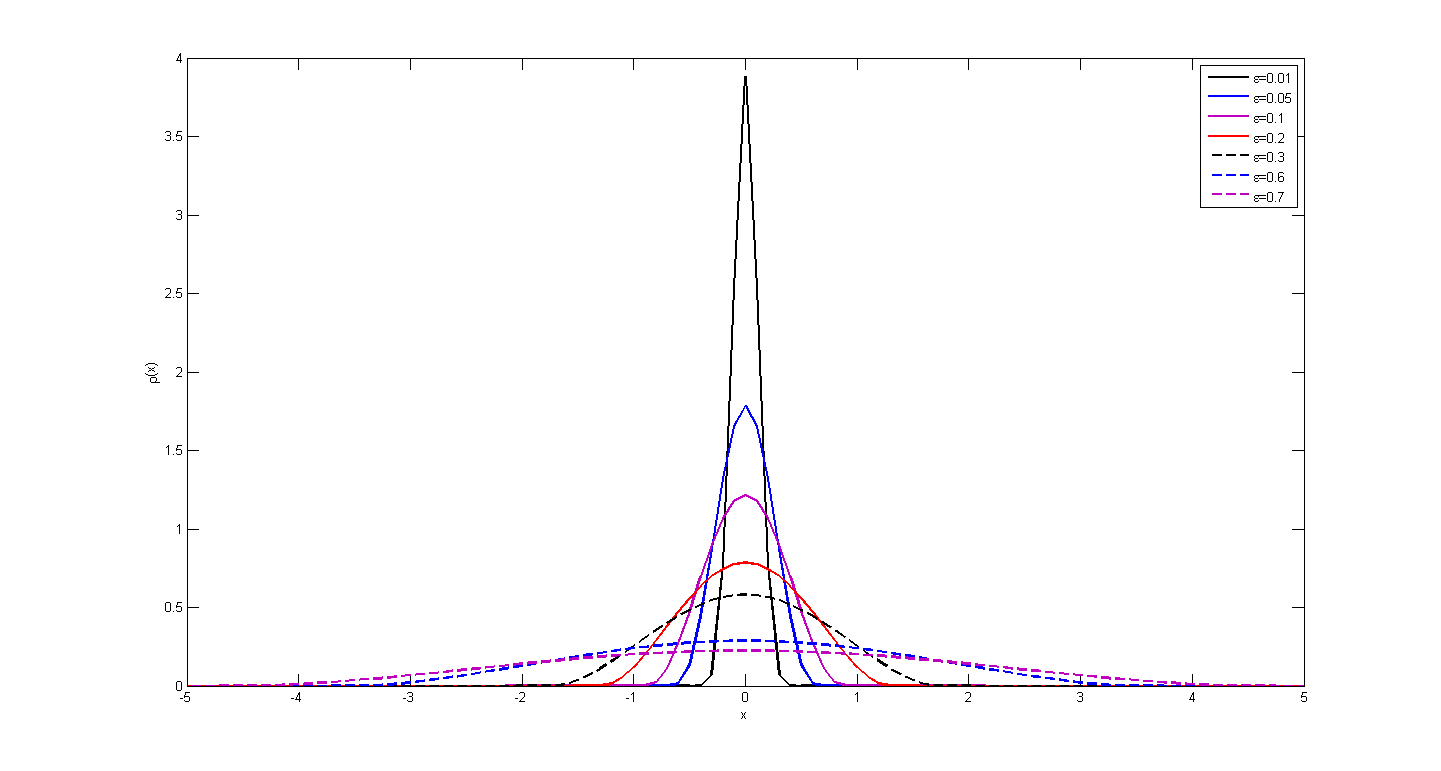}}\\
			\subfigure[]{ \includegraphics[width=0.35\textwidth]{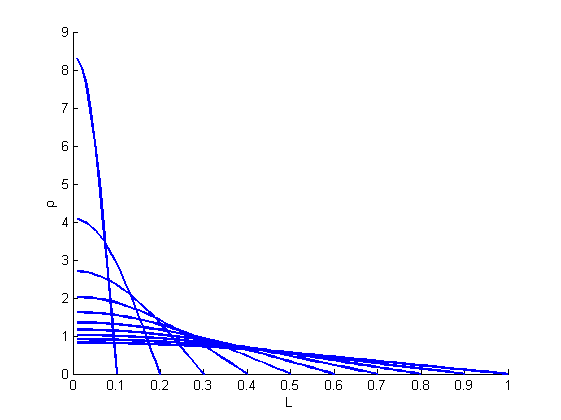}}
	\subfigure[]{	\includegraphics[width=0.35\textwidth]{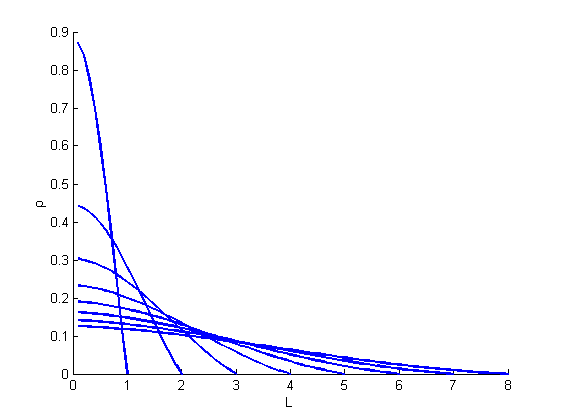}}
	\caption{(a) Kernel $G(x)= \frac{1}{2}\exp(-\left|x\right|)$; (b) Stationary solutions for equation \eqref{eq:main_evo}  for the kernel $G$.
	         (c),(d) eigenfunctions for $\epsilon(L)$ with $L\in [0,1]$ resp. $L\in [1,8]$.}
		\label{fig:eigenfunctionsKonvex}
\end{figure}

Furthermore we computed the stationary solutions \eqref{eq:main_evo} with the kernel $G(x)= \frac{1}{2}\exp(-\left|x\right|)$  (see Figure \ref{fig:eigenfunctionsKonvex}(a)) with Lipschitz singularity at the point zero. We present the results in Figure \ref{fig:eigenfunctionsKonvex}(b) for different $\epsilon$. Additionally we present eigenfunctions for the corresponding operator $\mathcal{G}_L$ \ref{eg:Gl} in Figure \ref{fig:eigenfunctionsKonvex}(c),(d).

\section*{Acknowledgement}

Most of the work of this paper was carried out during two visits by MDF at the Institute for Computational and Applied Mathematics of the University of Muenster. He is very grateful to MB and his group for their support and hospitality. MB and MDF thank the Isaac Newton Institute, Cambridge, for hospitality and financial support during the programme {\em Partial Differential Equations in Kinetic Theory}.

\bibliography{aggregation_diffusion}
\bibliographystyle{abbrv}

\end{document}